\documentclass[10pt]{article}
\usepackage{amssymb}
\usepackage{amsmath}

\usepackage[hypertex]{hyperref}
\usepackage{url}
\newtheorem{theo}{Theorem}[section]
\newtheorem{prop}[theo]{Proposition}
\newtheorem{lemm}[theo]{Lemma}
\newtheorem{coro}[theo]{Corollary}
\newtheorem{rema}[theo]{Remark}

\newtheorem{conj}[theo]{Conjecture}
\newtheorem{notation}[theo]{Notation}
\newtheorem{example}[theo]{Example}
\begin{document}
\author{}
\date{}

\newcommand{\cqfd}
{%
\mbox{}%
\nolinebreak%
\hfill%
\rule{2mm}{2mm}%
\medbreak%
\par%
}
\newfont{\gothic}{eufb10}

\title{The generalized Hodge  and Bloch conjectures  are equivalent for  general complete intersections}

\author{Claire Voisin
\\CNRS, Institut de math\'{e}matiques de Jussieu}

 \maketitle \setcounter{section}{-1}
 \begin{abstract} Let $X$ be a smooth complex projective variety
with trivial Chow groups. (By trivial, we mean that the cycle class is injective.)
We show (assuming the Lefschetz standard conjecture) that if the vanishing cohomology of  a general complete intersection $Y$  of ample hypersurfaces in $X$
has geometric coniveau $\geq c$, then the Chow groups  of cycles of dimension $\leq c-1$ of $Y$ are trivial.
The generalized Bloch conjecture for $Y$ is this statement
with ``geometric coniveau" replaced by ``Hodge coniveau".
\end{abstract}
\section{Introduction}
\setcounter{equation}{0}
  Recall first
 that a weight $k$  Hodge structure $(L,L^{p,q})$ has coniveau $c\leq\frac{k}{2}$ if the Hodge decomposition of $L_\mathbb{C}$ takes the form
 $$L_\mathbb{C}=L^{k-c,c}\oplus L^{k-c-1,c+1}\oplus\ldots\oplus L^{c,k-c}$$
 with $L^{k-c,c}\not=0$.
 If $X$ is a smooth complex projective variety and $Y\subset X$ is
 a closed algebraic subset of codimension $c$, then ${\rm Ker}\,(H^k(X,\mathbb{Q})\rightarrow H^k(X\setminus Y,\mathbb{Q}))$
 is a sub-Hodge structure of coniveau $\geq c$ of $H^k(X,\mathbb{Q})$ (cf. \cite[Theorem 7]{voisintorino}).
 The  generalized Hodge conjecture  formulated by Grothendieck \cite{grothendieck} is the following.
\begin{conj} \label{hodgeconjgen} Let $X$ be as above and let $L\subset H^k(X,\mathbb{Q})$ be a sub-Hodge structure of coniveau $\geq c$.
Then there exists  a closed algebraic subset $Y\subset X$ of codimension $c$ such that
$$L\subset {\rm Ker}\,(H^k(X,\mathbb{Q})\rightarrow H^k(X\setminus Y,\mathbb{Q})).$$
\end{conj}
This conjecture is widely open, even for general hypersurfaces or complete intersections in projective space
(cf. \cite{voisingafa}). Consider a smooth complete intersection $X\subset \mathbb{P}^n$ of $r$  hypersurfaces
of degrees $d_1\leq\ldots\leq d_r$. Then the coniveau of the Hodge structure on
$H^{n-r}(X,\mathbb{Q})_{prim}$ (the only part of the cohomology of $X$ which does not come from the ambient space)
is given by the formula (cf. \cite{voisingafa}, where  complete intersections of coniveau $2$ are studied):
\begin{eqnarray}\label{computconiveau}
{\rm coniveau}(H^{n-r}(X,\mathbb{Q})_{prim})\geq c\Leftrightarrow n\geq\sum_{i}d_i+(c-1)d_r.
\end{eqnarray}

The importance  of Conjecture \ref{hodgeconjgen}   has been underlined by the various  generalizations of  Mumford theorem
obtained in \cite{blochsrinivas}, \cite{schoen}, \cite{lewis}, \cite{laterveer}, \cite{paranjape}, based
on  refinements of the diagonal decomposition principle due to Bloch and Srinivas.
The resulting statement is the following (see \cite[II,10.3.2]{voisinbook}):
\begin{theo} \label{theotoutlemonde} Let $X$ be a smooth projective variety of dimension $m$. Assume that the cycle class map
$$cl: CH_i(X)_\mathbb{Q}\rightarrow H^{2m-2i}(X,\mathbb{Q})$$ is injective
for $i\leq c-1$. Then
we have $H^{p,q}(X)=0$ for $p\not=q$  and $p<c$ (or $q<c$). Hence the Hodge structures
on $H^{k}(X,\mathbb{Q})/H^{k}(X,\mathbb{Q})_{alg}$ are all of coniveau $\geq c$ and they satisfy the generalized Hodge conjecture
for coniveau $c$.
\end{theo}
Here $H^{k}(X,\mathbb{Q})_{alg}$ denotes the
``algebraic part'' of the cohomology, generated by classes of algebraic cycles. It is non zero of course only if
$k$ is even.

The first case of this theorem,  that is the case where $c=1$, was obtained by Bloch-Srinivas \cite{blochsrinivas}. It says that if a variety
$X$ has $CH_0(X)=\mathbb{Z}$, then $H^{k,0}(X)=0$ for any $k>0$ (which generalizes Mumford's theorem \cite{mumford})
and furthermore, the cohomology of positive degree of $X$ is supported on a proper algebraic subset $Y\subset X$ (which solves
Conjecture \ref{hodgeconjgen} for such $X$ and for coniveau $1$).

The next  major open problem, which by the above theorem would solve the generalized Hodge conjecture, is the following
conjecture relating the Hodge coniveau and Chow groups. This converse of Theorem \ref{theotoutlemonde}
is  a vast generalization of
Bloch conjecture for surfaces \cite{blochbook}.
\begin{conj} \label{blochconj} (cf. \cite[II, 11.2.2]{voisinbook}) Let $X$ be a smooth projective variety of dimension $m$ satisfying the condition
$H^{p,q}(X)=0$ for $p\not=q$ and $p<c$ (or $q<c$). Then for any
integer $i\leq c-1$, the cycle map
$$cl: CH_i(X)_\mathbb{Q}\rightarrow H^{2m-2i}(X,\mathbb{Q})$$ is injective.
\end{conj}
If we look at the case of hypersurfaces or complete intersections in projective space,
we see from (\ref{computconiveau}) that Conjecture \ref{blochconj} predicts the following :
\begin{conj}\label{conjhyp} Let $X\subset \mathbb{P}^n$ be a smooth complete intersection of hypersurfaces
of degrees $d_1\leq\ldots\leq d_r$. Then if
$n\geq\sum_{i}d_i+(c-1)d_r$,
the cycle map
$$cl: CH_i(X)_\mathbb{Q}\rightarrow H^{2n-2r-2i}(X,\mathbb{Q})$$ is injective for any
integer $i\leq c-1$.
\end{conj}
Note that, by Theorem \ref{theotoutlemonde}, this conjecture would imply Conjecture \ref{hodgeconjgen} for very general complete intersections.
Indeed, by monodromy arguments, the Hodge structure on the primitive middle cohomology of a very general
complete intersection is simple except for some rare and classified cases where it is made of Hodge classes. Thus
a sub-Hodge structure in this case must be the whole primitive cohomology in this case, and its coniveau is computed
by (\ref{computconiveau}).

Apart from very particular values of the degrees
$d_i$ (like complete intersections of quadrics \cite{otwi},
or cubics of small dimension \cite{collino}),
Conjecture \ref{conjhyp} is essentially  known only in the  case $c=1$, where the considered complete intersections are Fano, hence rationally connected,
so that the equality $CH_0(X)=\mathbb{Z}$ is trivial in this case.

In the paper \cite{voisincras}, it is proved that for any pair $(n,d)$, there are smooth hypersurfaces of degree
$d$ in $\mathbb{P}^n$ satisfying the conclusion of Conjecture \ref{conjhyp}.

 Coming back to Conjecture \ref{hodgeconjgen} for general complete intersections in projective space, we get from
 (\ref{computconiveau}) that it is equivalent  to the following statement:
 \begin{conj} \label{conjconiveauhyp}
 The primitive cohomology $H^{n-r}(X,\mathbb{Q})_{prim}$
 of a smooth complete intersection $X\subset \mathbb{P}^n$ of $r$  hypersurfaces of degrees $d_1\leq\ldots\leq d_r$
 vanishes on the complement of a closed algebraic subset
 $Y\subset X$ of codimension $c$ if
 $ n\geq\sum_{i}d_i+(c-1)d_r$.
 \end{conj}
 As already mentioned,  Conjecture \ref{conjconiveauhyp} would be implied by  Conjecture \ref{conjhyp} by Theorem \ref{theotoutlemonde}.
 The paper \cite{voisingafa} is an attempt to prove directly
 Conjecture \ref{conjconiveauhyp} for hypersurfaces or complete intersections of coniveau $2$ without trying to show the triviality of their
 $CH_0$ and $CH_1$ groups.
   Conjecture \ref{conjconiveauhyp}
  is shown there  to be implied by a conjecture on the effective cone of algebraic cycles (on some auxiliary variety).
 This work was motivated by the fact
that, unlike the case of coniveau $1$, and as is apparent from the lack of progresses in this direction  and the relative weakness
of the results obtained this way (see \cite{evl}, \cite{paranjape}, \cite{otw}), it seems now unlikely that  one will prove  Conjecture \ref{conjconiveauhyp}
by an application of Theorem
\ref{theotoutlemonde}, that is
 via the proof of the triviality of Chow groups of small dimension.

In fact, we will  essentially show in this paper that for general complete intersections inside any smooth projective variety
$X$ with ``trivial'' Chow groups,
Conjecture \ref{hodgeconjgen} (that is Conjecture
\ref{conjconiveauhyp} if $X=\mathbb{P}^n$)
 implies Conjecture \ref{blochconj} (that is Conjecture
\ref{conjhyp} if $X=\mathbb{P}^n$). Stated this way, this is not completely
correct, and we have to add an extra  assumption that we now explain.

Let us state the following conjecture, that we will relate later on (cf. Proposition
\ref{propcompconj}) to the so-called standard conjectures \cite{kleiman}:
\begin{conj}\label{conjstandard} Let $X$ be a smooth complex algebraic variety, and let $Y\subset X$ be
a closed algebraic subset. Let $Z\subset X$ be a codimension $k$ algebraic cycle, and assume that
the cohomology class $[Z]\in H^{2k}(X,\mathbb{Q})$ vanishes in $H^{2k}(X\setminus Y,\mathbb{Q})$. Then there exists
a codimension $k$ cycle $Z'$ on $X$ with $\mathbb{Q}$-coefficients, which is supported on $Y$ and such that $[Z']=[Z]$ in $H^{2k}(X,\mathbb{Q})$.
\end{conj}
Our main result in this paper is the following.
\begin{theo} \label{theobrutintro} Assume conjecture \ref{conjstandard} holds for cycles of codimension
$n-r$. Let $X$ be a smooth complex projective variety satisfying the property
that the cycle map $cl: CH^i(X)_\mathbb{Q}\rightarrow H^{2i}(X,\mathbb{Q})$ is injective for any $i$.
Let $L_1,\ldots, L_r,\,r\leq {\rm dim}\,X$, be  very ample line bundles on $X$. Assume that for a very general
complete intersection $X_b=X_1\cap\ldots\cap X_r$ of hypersurfaces $X_i\in | L_i|$, the
Hodge structure on $H^{n-r}(X_b,\mathbb{Q})_{prim}$ is supported on a closed algebraic subset $Y_b\subset X_b$ of codimension
$\geq c$.
Then for the general such  $X_b$ (hence in fact for all), the cycle map $cl:CH^i(X_b)\rightarrow H^{2i}(X_b,\mathbb{Q})$ is injective for any $i<c$.
\end{theo}
As a particular case, we get:
\begin{coro} \label{corointro} Assuming conjecture \ref{conjstandard} for cycles
of codimension $n-r$, the generalized Hodge conjecture \ref{conjconiveauhyp} implies the generalized Bloch conjecture
\ref{conjhyp}
for complete intersections in projective space.
\end{coro}

These results are conditional results. However
there is one non trivial case where they are unconditional, namely
the case of surfaces, where we get the following.
\begin{theo} \label{surfaces}  Let $X$ be a smooth complex projective variety of dimension $r+2$ satisfying the property
that the cycle map $cl: CH^i(X)_\mathbb{Q}\rightarrow H^{2i}(X,\mathbb{Q})$ is injective for any $i$.
Let $L_1,\ldots, L_r,\,r\leq {\rm dim}\,X$, be  very ample line bundles on $X$. Assume the smooth
complete intersections surfaces  $X_b=X_1\cap\ldots\cap X_r$, $X_i\in | L_i|$ have $h^{2,0}(X_b)=0$.
Then we have $CH_0(X_b)=\mathbb{Z}$.
\end{theo}
{\bf Proof.} Indeed, the assumption $h^{2,0}(X_b)=0$ implies by the Lefschetz theorem
on $(1,1)$-classes that
the Hodge structure on $H^2(X_b,\mathbb{Q})$ is generated by divisor classes.
So, the generalized Hodge conjecture is true in this case.
Furthermore, Conjecture \ref{conjstandard} will be  satisfied in this case, because
it is satisfied by codimension $2$ cycles (cf. Lemma \ref{lemmacodim2}). Finally, $CH_0(X_b)_\mathbb{Q}=\mathbb{Q}$
implies $CH_0(X_b)=\mathbb{Z}$ by Roitman's theorem \cite{roitman}.

\cqfd
We refer to Remark \ref{remasurfacesthreefolds} for the case of threefolds, where we get a similar conclusion.

The paper is organized as follows: in section \ref{sec1}, we will show that Conjecture \ref{conjstandard} is implied by the so-called Lefschetz
conjecture.
In section \ref{sec2}, we will prove Theorem \ref{theobrutintro}.
Section \ref{sec3} will provide a number of other geometric applications. For example, we will show how to recover
the results of \cite{voisinscuola}, or \cite{peters}.
 We will get  more generally results for many complete intersections
  $X_b$ endowed with the action of a finite group $G$. In this case, the method applies as well to the $\chi$-invariant part
 of $CH(X_b)$ where $\chi:G\rightarrow\{1,-1\}$ is a character.
 \begin{rema}{\rm As already mentioned, our  results are  unconditional in the case of surfaces, where the
 needed assumptions will be  satisfied by Lefschetz's theorem on $(1,1)$-classes. This is uninteresting
 in the  case of complete intersections surfaces in projective space, since those of coniveau $1$ are Fano, but in the presence of a group action,
 there are interesting non trivial examples of group actions on complete
 intersection surfaces  where this applies, in particular the Godeaux surfaces considered in \cite{voisinscuola}.}
 \end{rema}
Another potential application concerns self-products of Calabi-Yau hypersurfaces. It was noticed
in \cite{voisinsymmetric} that the generalized Bloch conjecture implies the following: Let
$X$ be a $n$-dimensional smooth projective variety with $H^i(X,\mathcal{O}_X)=0$ for $0<i<n$ and
$H^n(X,\mathcal{O}_X)=\mathbb{C}$ (for example $X$ could be a Calabi-Yau manifold).
Then if $n$ is even the antisymmetric product $z\times z'-z'\times z$ of two $0$-cycles of
$X$ of degree $0$ should be $0$  in $CH_0(X\times X)$. If $n$ is odd the symmetric product $z\times z'+z'\times z$ of two $0$-cycles of
$X$ of degree $0$ should be $0$  in $CH_0(X\times X)$. This comes from Conjecture
\ref{blochconj}, or rather its generalization to motives,
and from the observation
that the Hodge structure on $\bigwedge ^2 H^n(X,\mathbb{Q})$ has coniveau $\geq 1$ (see Lemma \ref{leconiveauCY}).

We will show in Section \ref{sec3}:
\begin{theo} \label{theoproduitCY} Assume Conjecture \ref{conjstandard}. Let $X$ be a   Calabi-Yau
hypersurface in projective space $\mathbb{P}^n$. Then if the generalized Hodge conjecture
is true for the coniveau $1$ Hodge structure on $\bigwedge^2H^{n-1}(X,\mathbb{Q})$
(seen as a sub-Hodge structure of $H^{2n-2}(X\times X,\mathbb{Q})$) for the general such $X$,
the antisymmetric product $z\times z'-z'\times z$ of two $0$-cycles of
$X$ of degree $0$ is equal to  $0$  in $CH_0(X\times X)$ for $n-1$ even and the symmetric product $z\times z'+z'\times z$ of two $0$-cycles of
$X$ of degree $0$ is equal to $0$  in $CH_0(X\times X)$ for $n-1$ odd.
\end{theo}
{\bf Thanks.}  I thank Manfred Lehn, Christoph Sorger  and Burt Totaro for their help concerning
Lemma \ref{ledesing}, and especially Manfred Lehn and  Christoph Sorger for communicating
the text \cite{lehnsorger}.
\section{Remarks on Conjecture \ref{conjstandard}\label{sec1}}
The aim of this section is to comment on
Conjecture \ref{conjstandard}. The first observation to make is the following:
\begin{lemm}\label{lemmacodim2}
 Conjecture \ref{conjstandard} is satisfied by codimension $k$ cycles  whose cohomology class vanishes
away from a codimension $k-1$ closed algebraic subset.

In particular, Conjecture \ref{conjstandard} is satisfied by codimension $2$ cycles.
\end{lemm}
{\bf Proof.}   Indeed, if we have a codimension $k$ cycle $Z\subset X$, whose cohomology class
$[Z]\in H^{2k}(X,\mathbb{Q})$ vanishes on the open set $X\setminus Y$, where ${\rm codim\,Y}\geq k-1$,
then
we know (cf. \cite[Proposition 3]{voisintorino}) that there are  Hodge classes $\alpha_i\in Hdg^{2k-2c_i}(\widetilde{Y}_i,\mathbb{Q})$, such that
$$[Z]=\sum_i\tilde{j}_{i*}\alpha_i,$$
where $\tilde{j}_{i}:\widetilde{Y}_i\rightarrow X$ are desingularizations of
the irreducible components $Y_i$ of $Y$, and $c_i:={\rm codim}\,Y_i$.

As $c_i\geq k-1$ for all $i$'s, the classes $\alpha_i$ are  cycle classes on $\widetilde{Y}_i$
by the Lefschetz theorem on $(1,1)$-classes, which concludes the proof.
\cqfd

We are now going to   relate precisely Conjecture \ref{conjstandard} to the famous ``standard conjectures'' \cite{kleiman}.
Let $X$ be a smooth projective variety of dimension $n$. The K\"{u}nneth decomposition of $H^*(X\times X,\mathbb{Q})$ gives:
$$ H^m(X\times X,\mathbb{Q})\cong \bigoplus_{p+q=m}H^p(X,\mathbb{Q})\otimes H^q(X,\mathbb{Q}).$$
 Poincar\'{e} duality on $X$ allows to rewrite this as
\begin{eqnarray}
\label{kundec}
H^m(X\times X,\mathbb{Q})\cong \bigoplus_{p+q=m}{\rm Hom}\,(H^{2n-p}(X,\mathbb{Q}), H^q(X,\mathbb{Q})).
\end{eqnarray}
On the other hand, we have the following lemma (cf. \cite[I, 11.3.3]{voisinbook}):
\begin{lemm} \label{lehoco}
Let $m=p+q$ be even. A cohomology class $$\alpha\in {\rm Hom}\,(H^{2n-p}(X,\mathbb{Q}), H^q(X,\mathbb{Q}))\subset H^m(X\times X,\mathbb{Q})$$ is
a Hodge class on $X\times X$ if and only if it is a morphism of Hodge structures.
\end{lemm}
There are two kinds of particularly interesting Hodge classes on $X\times X$ obtained from Lemma \ref{lehoco}.

a)  Let $m=2n$ and consider for each $0\leq q\leq 2n$ the element
$ Id_{H^q(X,\mathbb{Q})}$ which provides by (\ref{kundec}) and Lemma \ref{lehoco} a Hodge class $\delta_q$ of degree $2n$ on
$X$. This class is called the $q$-th K\"{u}nneth component of the diagonal of $X$.
The first standard conjecture (or K\"{u}nneth's conjecture,  cf.  \cite{kleiman})
is the following:
\begin{conj}\label{kunconj}  The classes $\delta_i$ are algebraic, that is, are classes of algebraic cycles on $X\times X$ with rational
coefficients.
\end{conj}
b) Let $L$ be an ample line bundle on $X$, and $l:=c_1(L)\in H^2(X,\mathbb{Q})$.
For any integer $k\leq n$, the hard Lefschetz theorem \cite[I, 6.2.3]{voisinbook} says that
the cup-product map
$$ l^{n-k}\cup: H^k(X,\mathbb{Q})\rightarrow H^{2n-k}(X,\mathbb{Q})$$
is an isomorphism. This is clearly an isomorphism of Hodge structure. Its inverse
$$(l^{n-k}\cup)^{-1}:H^{2n-k}(X,\mathbb{Q})\rightarrow H^k(X,\mathbb{Q})$$
is also an isomorphism of Hodge structures, which by (\ref{kundec}) and Lemma \ref{lehoco}
provides a Hodge class $\lambda_{n-k}$ of degree $2k$ on $X\times X$.
The second standard conjecture we will consider (this is one form of Lefschetz' conjecture, cf.  \cite{kleiman})
is the following:
\begin{conj}\label{lefconj}  The classes $\lambda_i$ are algebraic, that is, are classes of algebraic cycles on $X\times X$ with rational
coefficients.
\end{conj}
\begin{rema}{\rm One could also ask whether there is a codimension $n$ algebraic cycle
$Z$ on $X\times X$ with rational
coefficients such that the induced morphism $[Z]_*:H^{2n-k}(X,\mathbb{Q})\rightarrow H^k(X,\mathbb{Q})$
is equal to $\lambda_k$. However, if this weaker version is true for any $k$, the K\"{u}nneth decomposition
is algebraic and then by taking the K\"{u}nneth component of bidegree $(k,k)$ of $Z$, we get  an affirmative answer
to Conjecture \ref{lefconj}.}
\end{rema}

Let us show the following
\begin{prop} \label{propcompconj} The Lefschetz conjecture for any $X$ is equivalent to the conjunction of the K\"{u}nneth conjecture
\ref{kunconj} and of
Conjecture \ref{conjstandard} for any $X$.
\end{prop}
{\bf Proof.} Let us assume that the K\"{u}nneth conjecture holds for  $X$ and Conjecture \ref{conjstandard} holds for any
pair $Y\subset X'$. Let $i<n$. Consider the K\"{u}nneth component $\delta_{2n-i}$ of $\Delta_X$, so
$\delta_{2n-i}\in H^{i}(X,\mathbb{Q})\otimes H^{2n-i}(X,\mathbb{Q})$ is the class of an algebraic cycle $Z$
on $X\times X$.
Let $Y_i\stackrel{j_i}{\hookrightarrow} X$ be a smooth complete intersection of $n-i$ ample hypersurfaces in $X$. Then
the Lefschetz theorem on hyperplane sections (cf. \cite[II, 1.2.2]{voisinbook})
says that
$$j_{i*}: H^{i}(Y_i,\mathbb{Q})\rightarrow H^{2n-i}(X,\mathbb{Q})$$
is surjective.
It follows that the class of the cycle $Z$ vanishes on $X\times (X\setminus Y_i)$.
By Conjecture \ref{conjstandard}, there is a $n$-cycle $Z'$ supported on
$X\times Y_i$ such that the class  $(id,j)_*[Z']$ is equal to $[Z]$.
Consider the morphism of Hodge structures induced by $[Z']$:
$$[Z']_*: H^{2n-i}(X,\mathbb{Q})\rightarrow H^i(Y_i,\mathbb{Q}).$$
Composing with the morphism $j_{i*}: H^i(Y_i,\mathbb{Q})\rightarrow H^{2n-i}(X,\mathbb{Q})$,
we get $j_{i*}\circ[Z']_*=Id_{H^{2n-i}(X,\mathbb{Q})}$. It follows that $[Z']_*$ is injective, and that its transpose
$[Z']^*: H^i(Y_i,\mathbb{Q})\rightarrow H^i(X,\mathbb{Q})$ is surjective.
We now apply \cite[Proposition 8]{charles} and induction on dimension to conclude that Lefschetz' conjecture holds for  $X$.

Conversely, assume Lefschetz' conjecture holds for any smooth projective variety. It obviously implies the
K\"{u}nneth conjecture. Let  now  $X, \,Y,\,Z$ be as in Conjecture \ref{conjstandard}. Set
$n={\rm dim}\,X,\,k={\rm codim}\,Z,\,l={\rm codim}\,Y$.
Let $\tilde{j}:\widetilde{Y}\rightarrow X$ be a desingularization
of $Y$. It is known (cf. \cite[Proposition 3]{voisintorino}) that there exists a Hodge class $\beta\in H^{2k-2l}(\widetilde{Y},\mathbb{Q})$
such that $\tilde{j}_*\beta=[Z]$. The question is to find such a $\beta$ algebraic on
$\widetilde{Y}$. The argument is easier to understand if we assume that ${\rm dim}\, X=2{\rm dim}\,Z$ and ${\rm dim}\,\widetilde{Y}={\rm dim}\,X$,
which can always be achieved up to replacing $(X,\,Z)$ by $(X\times \mathbb{P}^{l_1}, Z\times \mathbb{P}^{l_2})$ with
${\rm dim}\,X+l_1=2({\rm dim}\,Z+l_2)$, and
$Y$ by $Y\times \mathbb{P}^{l_3}$, with ${\rm dim}\,Y+l_3={\rm dim}\,X+l_1$.
Let thus $2n={\rm dim}\,X={\rm dim}\,\widetilde{Y}$.
The Lefschetz conjecture for $\widetilde{Y}$ implies that the Lefschetz decomposition
$$H^{2n}(\widetilde{Y},\mathbb{Q})=\oplus_{n-k\geq0} h_{\widetilde{Y}}^k \cup H^{2n-2k}(\widetilde{Y},\mathbb{Q})_{prim},$$
(where $h_{\widetilde{Y}}\in H^2(\widetilde{Y},\mathbb{Q})$ is the class of an ample hypersurface on $\widetilde{Y}$,)
is algebraic  in the sense that the various projectors $\pi_k$ on primitive factors are induced by
classes of algebraic cycles on $\widetilde{Y}\times\widetilde{Y}$.
Let $G:=\sum_k(-1)^k\pi_k$. $G$ is thus induced by an algebraic cycle of dimension
$2n$ on $\widetilde{Y}\times\widetilde{Y}$. The Hodge-Riemann bilinear relations (\cite[I, 6.3.2]{voisinbook})
say that the symmetric intersection pairing
$$<\alpha,\beta>_G=<\alpha,G_*\beta>$$
is nondegenerate of a definite sign on the space $H^{n,n}(\widetilde{Y})\cap H^{2n}(\widetilde{Y},\mathbb{R})$
of real cohomology classes of Hodge type $(n,n)$
 on $\widetilde{Y}$, hence in particular on
the space $Hdg^{2n}(\widetilde{Y})$ of Hodge classes on $\widetilde{Y}$.
 It follows that
 the image of
 $$\widetilde{j}_*: Hdg^{2n}(\widetilde{Y},\mathbb{Q})\rightarrow Hdg^{2n}(X,\mathbb{Q})$$
 is the same as the image of
 $$ \widetilde{j}_*\circ G\circ \widetilde{j}^*: Hdg^{2n}(X,\mathbb{Q})\rightarrow Hdg^{2n}(X,\mathbb{Q}).$$
It follows that $\widetilde{j}_*\circ G\circ \widetilde{j}^*$ induces an automorphism  of
the subspace $\widetilde{j}_*Hdg^{2n}(\widetilde{Y},\mathbb{Q})\subset Hdg^{2n}(X,\mathbb{Q})$.
This automorphism is induced by an algebraic self-correspondence of $X$, that is a $2n$-algebraic cycle of $X\times X$
with $\mathbb{Q}$-coefficients. Hence it preserves the subspace of algebraic
classes
$$(\widetilde{j}_*Hdg^{2n}(\widetilde{Y},\mathbb{Q}))\cap H^{2n}(X,\mathbb{Q})_{alg}\subset Hdg^{2n}(X,\mathbb{Q}).$$
 It thus also induces an automorphism of this subspace.

Recalling that $[Z]=\widetilde{j}_*\beta\in {\rm Im}\,\widetilde{j}_*\in \widetilde{j}_*Hdg^{2n}(\widetilde{Y},\mathbb{Q})\cap H^{2n}(X,\mathbb{Q})_{alg}$, we thus conclude that
$[Z]=\widetilde{j}_*((G\circ \widetilde{j})^*\gamma)$
for some $\gamma\in H^{2n}(X,\mathbb{Q})_{alg}$.
This concludes the proof since
 $(G\circ\widetilde{j}^*) \gamma$ is an algebraic class on $\widetilde{Y}$.
\cqfd
The main use we will make of Conjecture \ref{conjstandard} is the following strengthening of the generalized Hodge conjecture.
Let $X$ be a smooth complex projective variety of dimension $n$, and let
$L$ be a sub-Hodge structure of $H^n(X,\mathbb{Q})_{prim}$, where the subscript ``prim'' stands for ``primitive with respect to a given
polarization on $X$''. We know then by the second Hodge-Riemann bilinear relations
\cite[I, 6.3.2]{voisinbook} that the intersection form $<,>$ restricted to $L$ is nondegenerate. Let $\pi_L:H^n(X,\mathbb{Q})\rightarrow L$
be the orthogonal projector on $L$. We assume that $\pi_L$ is algebraic, that is,  there is
a  $n$-cycle $\Delta_L\subset X\times X$, such that
$$[\Delta_L]_*=\pi_L:H^n(X,\mathbb{Q})\rightarrow L\subset H^n(X,\mathbb{Q}),$$
$$[\Delta_L]_*=0:H^i(X,\mathbb{Q})\rightarrow  H^i(X,\mathbb{Q}),\,\,i\not=n.$$
\begin{lemm} \label{leeasybutcrucial} Assume that there exists a closed algebraic subset $Y\subset X$ such that
$L$ vanishes in $H^n(X\setminus Y,\mathbb{Q})$. Then if Conjecture \ref{conjstandard} holds,
there is a cycle $Z'_L\subset Y\times Y$ with $\mathbb{Q}$-coefficients  such that
$$[Z'_L]=[\Delta_L]\,\,{\rm in}\,\,H^{2n}(X\times X,\mathbb{Q}).$$

\end{lemm}

{\bf Proof.} Indeed,  because $\pi_L$ is the orthogonal projector on $L$, the class $[Z_L]$ belongs to
$L\otimes L\subset H^{2n}(X\times X,\mathbb{Q})$. As $L$ vanishes in $H^n(X\setminus Y,\mathbb{Q})$,
the class $[Z_L]\in L\otimes L$ vanishes in $H^{2n}(X\times X\setminus (Y\times Y),\mathbb{Q})$.
Conjecture \ref{conjstandard} then guarantees the existence of a cycle
$Z'_L\subset Y\times Y$ such that
$[Z'_L]=[\Delta_L]\,\,{\rm in}\,\,H^{2n}(X\times X,\mathbb{Q})$.

\cqfd
\section{Proof of Theorem \ref{theobrutintro}\label{sec2}}
\subsection{Generalities on varieties with ``trivial'' Chow groups}
We will say that a (non necessarily projective) smooth variety  satisfies property $\mathcal{P}$ (or has trivial Chow groups)
if
the cycle map
$$cl:CH^i(X)_\mathbb{Q}\rightarrow H^{2i}(X,\mathbb{Q})$$
is injective for all $ i$.
We will say that $X$ satisfies property $\mathcal{P}_k$ if the cycle class map
$$cl:CH^i(X)_\mathbb{Q}\rightarrow H^{2i}(X,\mathbb{Q})$$
is injective for all $ i\leq k$.
When $X$ is projective, it is known (cf. \cite{laterveer}, \cite{lewis}) that if $X$ has trivial Chow groups,
the cycle class map
$$cl:CH^k(X)_\mathbb{Q}\rightarrow H^{2k}(X,\mathbb{Q})$$
is an isomorphism for any $k$ and that $H^{2k+1}(X,\mathbb{Q})=0$ for all $k$.
We have the following lemma:
\begin{lemm} \label{leU}Assume Conjecture \ref{conjstandard}. Let $X$ be a smooth projective
variety satisfying property $\mathcal{P}$. Then any Zariski open set
$U\subset X$ satisfies property $\mathcal{P}$.
\end{lemm}
{\bf Proof.} Write $U=X\setminus Y$. Let $Z$ be a codimension $k$ cycle on $U$ with vanishing cohomology class. Then
$Z$ is the restriction to $U$ of a cycle $\overline{Z}$ on $X$, which has the property that
$$[\overline{Z}]_{\mid U}=0.$$
Conjecture \ref{conjstandard} says that
there is a cycle $Z'$ supported on $Y$ such that $[\overline{Z}]=[Z']$ in $H^{2k}(X,\mathbb{Q})$.
The cycle $\overline{Z}-Z'$ is thus cohomologous to $0$ on $X$. As $X$ satisfies property $\mathcal{P}$,
$\overline{Z}-Z'$ is rationally equivalent to $0$ on $X$ modulo torsion, and so is its restriction
to $U$, which is equal to $Z$.
\cqfd
\begin{lemm} \label{leproj}Let $X$ be a smooth complex variety satisfying property $\mathcal{P}_k$. Then
any projective bundle $p:\mathbb{P}(E)\rightarrow X$, where $E$ is a locally free sheaf on $X$, satisfies property $\mathcal{P}_k$.
\end{lemm}
{\bf Proof.}  Indeed, let $h=c_1(\mathcal{O}_{\mathbb{P}(E)}(1))\in CH^1(\mathbb{P}(E))$
and let $[h]\in H^2(\mathbb{P}(E),\mathbb{Q})$ be its topological first Chern class.
The canonical decompositions (\cite[I,7.3.3]{voisinbook}, \cite[II,9.3.2]{voisinbook}
 $$CH^*(\mathbb{P}(E))_\mathbb{Q}=\oplus_{0\leq i\leq r-1} h^ip^* CH^{*-i}(X,\mathbb{Q}),$$
 $$H^*(\mathbb{P}(E))_\mathbb{Q}=\oplus_{0\leq i\leq r-1} [h]^i\cup p^* H^{*-2i}(X,\mathbb{Q}),$$
 are compatible with the cycle map $cl:CH^*(X)\rightarrow H^{2*}(X,\mathbb{Q})$.
 Thus if $cl$ is injective on cycles of codimension $\leq k$ on  $X$, it is also injective
on cycles of codimension $\leq k$ on   $\mathbb{P}(E)$.
\cqfd

We prove similarly.
 \begin{lemm} \label{leblowup} Let $X$ be a smooth complex algebraic variety satisfying property $\mathcal{P}_k$
 and let $Y\subset X$ be a subvariety satisfying
 property $\mathcal{P}_{k-1}$. Then
the blow-up $\widetilde{X}_Y\rightarrow  X$ of $X$ along $Y$
satisfies property $\mathcal{P}_k$.
\end{lemm}
Let us conclude with two more properties:
\begin{lemm} \label{ledominant} Let $\phi: X\rightarrow X'$ be a projective surjective morphism, where
$X$ and $X'$ are smooth complex algebraic variety. If $X$ satisfies property $\mathcal{P}$, so does $X'$.
\end{lemm}
{\bf Proof.} Indeed, let $h\in CH^1(X)$ be the first Chern class of a relatively ample line bundle. Let $r={\rm dim}\,X-{\rm dim}\, X'$, and let $d$ be defined
by $\phi_*h^r=d X'\in CH^0(X')$. Then we have the projection formula:
\begin{eqnarray}\label{projform}\phi_*(h^r\cdot\phi^*\alpha)=d\alpha,\,\forall \alpha\in CH^*(X')_\mathbb{Q}.
\end{eqnarray}
If $\alpha\in CH^*(X')_\mathbb{Q}$ satisfies $cl(\alpha)=0$ then $\phi^*(cl(\alpha))=cl(\phi^*\alpha)=0$
in $H^{2*}(X,\mathbb{Q})$. Thus $\phi^*\alpha=0$ in $CH^*(X)_\mathbb{Q}$ and $\alpha=0$
in $CH(X')_\mathbb{Q}$ by (\ref{projform}).
\cqfd

\begin{prop}\label{propproduit} Let $X$ be a smooth projective variety satisfying property $\mathcal{P}$. Then $X\times X$
satisfies property $\mathcal{P}$.
\end{prop}
{\bf Proof.} This uses the fact (proved eg in \cite{paranjape}) that a variety satisfying property $\mathcal{P}$ has
a complete decomposition of the diagonal as a combination of products of algebraic cycles (cf. also \cite[II,10.3.1]{voisinbook}):
$$\Delta_X=\sum_{i,j}n_{ij}Z_i\times Z_j\,\,{\rm in}\,\, CH^n(X\times X),$$
where $n_{ij}\in \mathbb{Q}$, and ${\rm dim}\,Z_i+{\rm dim}\,Z_j=n={\rm dim}\,X$.
It follows that the variety $Z:=X\times X$ also admits such decomposition,
since $\Delta_Z=p_{13}^*\Delta_X\cdot p_{24}^*\Delta_X$ in $CH^{2n}(Z\times Z)$, where
$p_{ij}$ is the projection of $Z\times Z=X^4$ to the product $X\times X$ of the $i$-th and  $j$-th summand.

But this in turn implies that $CH^*(Z)_\mathbb{Q}\cong H^{2*}(Z,\mathbb{Q})$. Indeed, write
$$\Delta_Z=\sum_{i,j}m_{ij} W_i\times W_j\,\,{\rm in}\,\, CH^{2n}(Z\times Z).$$
Then any cycle
$\gamma\in CH(Z)_\mathbb{Q}$ satisfies
$$\gamma=\Delta_{Z*}\gamma=\sum_{i,j}m_{ij} {\rm deg}\, (\gamma\cdot W_i)W_j\,\,{\rm in}\,\,CH(Z)_\mathbb{Q}.$$
It immediately follows that if $\gamma$ is homologous to $0$, it vanishes in $CH(Z)_\mathbb{Q}$.

\cqfd
\subsection{Proof of Theorem \ref{theobrutintro}}
We will start with a few preparatory Lemmas.
 Consider a smooth projective variety $X$ of dimension $n$ satisfying property $\mathcal{P}$. Let $L_i,\,i=1,\ldots,r$, be very ample line bundles on
$X$. Let $j:X_{b}\hookrightarrow X$ be a very general complete intersection of hypersurfaces in $|L_i|,\,i=1,\ldots,r$. Then
$X_b$ is smooth of dimension $n-r$, and the Lefschetz theorem on hyperplane sections says that
\begin{eqnarray}\label{decompvan} H^*(X_b,\mathbb{Q})=H^{n-r}(X_b,\mathbb{Q})_{van}\oplus H^*(X,\mathbb{Q})_{\mid X},
\end{eqnarray}
where the vanishing cohomology $L_b:=H^{n-r}(X_b,\mathbb{Q})_{van}$ is defined as
${\rm Ker}\,(j_*:H^{n-r}(X_b,\mathbb{Q})\rightarrow H^{n+r}(X,\mathbb{Q}))$. Note that
$H^{n-r}(X_b,\mathbb{Q})_{van}$ is contained in $H^{n-r}(X_b,\mathbb{Q})_{prim}$, where ``prim'' means primitive for a very ample line bundle
coming from $X$, and thus, by the second Hodge-Riemann
bilinear relations, the intersection form $<,>$ on $H^{n-r}(X_b,\mathbb{Q})$ remains nondegenerate after restriction to
$L_b$.

As $X$ has trivial Chow groups, we know that $H^*(X,\mathbb{Q})$ is generated by
classes of algebraic cycles and so is the restriction $H^*(X,\mathbb{Q})_{\mid X_b}$. This implies the following:
\begin{lemm}\label{lealgebraicpil}
 The orthogonal
projector $\pi_{L_b}$ on $L_b$ is algebraic.
\end{lemm}
{\bf Proof.} In fact, we can construct a canonical algebraic cycle
$\Delta_{b,van}$ with $\mathbb{Q}$-coefficients
on $X_b\times X_b$ whose class $[\Delta_{b,van}]$ is equal to
$\pi_{L_b}$. For this, we choose a basis of $\oplus_{i\leq n-r}H^{2i}(X,\mathbb{Q})$. As
we know that $X$ has trivial Chow groups, this basis consists of classes $[z_{i,j}]$ of algebraic cycles $z_{i,j}$
on $X$, with ${\rm codim}\,z_{i,j}=i\leq n-r$.
We choose an ample line bundle $H$ on $X$, and note that by the hard Lefschetz theorem, the classes
$[h]^{n-r-i}\cup [z_{i,j}]_{\mid X_b}$, together with the classes $[z_{i,j}]_{\mid X_b}$,
form a basis of $H^*(X,\mathbb{Q})_{\mid X_b}$. Here $[h]\in H^2(X_b,\mathbb{Q})$ is the topological first Chern class of $H_{\mid X_b}$.
The intersection form on $H^*(X_b,\mathbb{Q})$ is nondegenerate when restricted to
$H^*(X,\mathbb{Q})_{\mid X_b}$, and $L_{b}$ is the orthogonal complement
of $H^*(X,\mathbb{Q})_{\mid X_b}$ with respect to the intersection pairing
on $H^*(X_b,\mathbb{Q})$. We thus have the equality of
orthogonal projectors:
$$\pi_{L_b}+\pi_{H^*(X,\mathbb{Q})_{\mid X_b}}=Id_{H^*(X_b,\mathbb{Q})}.$$
But it is clear that the orthogonal projector $\pi_{H^*(X,\mathbb{Q})_{\mid X_b}}$ is given by the class of an algebraic cycle on
$X_b\times X_b$ (it suffices to choose an orthogonal basis of $H^*(X,\mathbb{Q})_{\mid X_b}$ for the intersection form
on $H^*(X_b,\mathbb{Q})$ and to recall that $H^*(X,\mathbb{Q})_{\mid X_b}$ consists of algebraic classes).
As $Id_{H^*(X_b,\mathbb{Q})}$ corresponds to the class of the diagonal of $X_b$, the proof is finished.
\cqfd

We now assume that there is a closed algebraic subset $Y_b\subset X_b$ of codimension $c$ such that $L_b$ vanishes on $X_b\setminus Y_b$.
Then, under Conjecture \ref{conjstandard}, Lemma \ref{leeasybutcrucial} tells that there is an algebraic cycle
$Z_b$ supported on $Y_b\times Y_b$ such that $[Z_b]=[\Delta_{b,van}]$.

The key point now is the following easy Proposition \ref{propcrux}. In the following, we will put everything in  family.
Let $\pi:\mathcal{X}\rightarrow B$ be a smooth projective morphism and let $(\pi,\pi):\mathcal{X}\times_B\mathcal{X}\rightarrow B$ be the
 fibered self-product of $\mathcal{X}$ over $B$.
Let $\mathcal{Z}\subset \mathcal{X}\times_B\mathcal{X}$ be a codimension $k$ algebraic cycle. We denote the fibres
$\mathcal{X}_b:=\pi^{-1}(b)$, $\mathcal{Z}_b:=\mathcal{Z}_{\mid \mathcal{X}_b\times \mathcal{X}_b}$.
\begin{prop}\label{propcrux}
Assume that for a very general point $b\in B$, there exists a closed algebraic subset
$Y_b\subset \mathcal{X}_b\times \mathcal{X}_b$ of codimension $c$, and
an algebraic cycle $Z'_b\subset Y_b\times Y_b$ with $\mathbb{Q}$-coefficients, such that
$$[Z'_b]=[\mathcal{Z}_b]\,\,{\rm in}\,\,H^{2k}(\mathcal{X}_b\times \mathcal{X}_b,\mathbb{Q}).$$
Then there exists a closed algebraic subset
$\mathcal{Y}\subset \mathcal{X}$ of codimension $c$, and a codimension $k$
 algebraic cycle $\mathcal{Z}'$ with $\mathbb{Q}$-coefficients on $\mathcal{X}\times_B\mathcal{X}$, which is supported
 on $\mathcal{Y} \times_B\mathcal{Y}$ and
such that for any $b\in B$,
$$[\mathcal{Z}'_b]=[\mathcal{Z}_b]\,\,{\rm in}\,\,H^{2k}(\mathcal{X}_b\times\mathcal{X}_b,\mathbb{Q}).$$
\end{prop}
\begin{rema}{\rm This proposition is a crucial observation in the present paper.  The key point
is  the fact that we do not need to make this base change
for this specific problem. The idea of spreading-out cycles has become very important in the theory of algebraic cycles
since Nori's paper \cite{nori}, (cf. \cite{greengriffiths}, \cite{ssaito}). For most problems however,
we usually need to work over a generically finite extension of the base, due to the fact that
cycles existing at the general point will exist on the total space of the family only
after a base change.

}
\end{rema}
{\bf Proof of Proposition \ref{propcrux}.} There are countably many algebraic varieties $M_i\rightarrow B$ parameterizing
data $(b,Y_b,Z'_b)$ as above, and we can assume that each $M_i$ parameterizes  universal objects
\begin{eqnarray}\label{families}
\mathcal{Y}_i\rightarrow M_i,\,\,\mathcal{Y}_i\subset \mathcal{X}_{M_i},\,\,
\mathcal{Z}'_i\subset \mathcal{Y}_i\times_{M_i}\mathcal{Y}_i,
\end{eqnarray}
satisfying the property that for $m\in M_i$, with $pr_1(m)=b\in B$,
$$ [\mathcal{Z}'_{i,b}]=[\mathcal{Z}_{i,b}]\,\,{\rm in}\,\,H^{2k}(\mathcal{X}_b\times \mathcal{X}_b,\mathbb{Q}).$$
 By assumption, $B$ is the union of the images of the first projections $M_i\rightarrow B$. By a Baire category
 argument, we conclude that one of the morphisms $M_i\rightarrow B$ is dominating.
 Taking a subvariety of $M_i$ if necessary, we may assume that
 $\phi_i:M_i\rightarrow B$ is generically finite. We may also assume that it is proper and carries the families
 $\mathcal{Y}_i\rightarrow M_i,\,\,\mathcal{Y}_i\subset \mathcal{X}_{M_i},\,\,
\mathcal{Z}'_i\subset \mathcal{Y}_i\times_{M_i}\mathcal{Y}_i$.
Denote by $r_i: \mathcal{X}_{M_i}\rightarrow \mathcal{X}$ the proper generically finite morphism
induced by $\phi_i$. Let
$$\mathcal{Y}:=r_i(\mathcal{Y}_i)\subset \mathcal{X}.$$
Note that because $r_i$ is generically finite, ${\rm codim}\,\mathcal{Y}\leq c$.
Let $r'_i:\mathcal{Y}_i\rightarrow \mathcal{Y}$ be the restriction of $r_i$ to $\mathcal{Y}_i$.
and let $\mathcal{Z}':=(r'_i,r'_i)_*(\mathcal{Z}'_i)$, which is a codimension $k$ cycle
in $\mathcal{X}\times_B\mathcal{X}$ supported in
$$(r'_i,r'_i)(\mathcal{Y}_i\times_{M_i}\mathcal{Y}_i)\subset \mathcal{Y}\times_B\mathcal{Y}.$$

 It is obvious that for any $b\in B$, $[\mathcal{Z}'_b]=N[\mathcal{Z}_b]\,\,{\rm in}\,\,H^{2k}(\mathcal{X}_b\times\mathcal{X}_b,\mathbb{Q})$,
 where $N$ is the degree of $r_i$.

\cqfd

In the application, $\mathcal{X}$ and $B$ will be the following :

\begin{notation}\label{notation}{\rm Let $X$ be a smooth projective with trivial Chow groups.
Let  $\mathbb{P}_i:=\mathbb{P}(H^0(X,L_i))$. Let
$B\subset\prod_i\mathbb{P}_i$ be the open set parameterizing smooth complete intersections and
let
$$\mathcal{X}\subset B\times X,\,\pi:\mathcal{X}\rightarrow B,$$
be the universal family. We will denote
${X}_b\subset \mathcal{X}$ the fibre $\pi^{-1}(b)$ for $b\in B$. }
 \end{notation}

We apply the previous proposition to
 $\mathcal{Z}=\mathcal{D}_{van}$,
the corrected relative diagonal with fibre  over $b\in B$ the $\Delta_{b,van}$ introduced in Lemma \ref{lealgebraicpil}. (Note that $\mathcal{D}_{van}$
is not in fact canonically defined, as it may be modified by adding cycles which are restrictions to
$\mathcal{X}$ of cycles in $CH^{>0}(B)\otimes CH(X)\subset CH(B\times X)$.)

We then get the following :
\begin{lemm} \label{coropropcrux} Assume that for a general point $b\in B$, there is a codimension
$c$ closed algebraic subset $Y_b\subset X_b$ such that $H^{n-r}(X_b,\mathbb{Q})_{van}$
vanishes on $X_b\setminus Y_b$. If furthermore Conjecture \ref{conjstandard} holds, there exists a closed algebraic subset
$\mathcal{Y}\subset \mathcal{X}$ of codimension $c$, and a codimension $n-r$
 algebraic cycle $\mathcal{Z}'$ on $\mathcal{X}\times_B\mathcal{X}$ with $\mathbb{Q}$-coefficients,
 which is supported on $\mathcal{Y}\times_B\mathcal{Y}$ and
such that for any $b\in B$,
$$[\mathcal{Z}'_b]=[\Delta_{b,van}]\,\,{\rm in}\,\,H^{2k}(X_b\times X_b,\mathbb{Q}).$$
\end{lemm}
{\bf Proof.} This is a direct application  of Proposition \ref{propcrux}, because we
know from Lemma \ref{leeasybutcrucial} that under Conjecture \ref{conjstandard}, the assumption implies that
there exists for any $b\in B$ an algebraic cycle $Z'_b\subset Y_b\times Y_b$ such that
$[Z'_b]=[\Delta_{b,van}]\,\,{\rm in}\,\,H^{2k}(X_b\times X_b,\mathbb{Q})$.
\cqfd

We have next the following :
\begin{lemm} \label{leleray} With notation as in  \ref{notation},
let $\alpha \in H^{2n-2r}(\mathcal{X}\times_B\mathcal{X},\mathbb{Q})$
be a cohomology class whose restriction to the fibres $X_b\times X_b$ is $0$.
Then we can write
$$\alpha=\alpha_1+\alpha_2$$
where $\alpha_1$ is the restriction to $\mathcal{X}\times_B\mathcal{X}$ of a
class $\beta_1\in H^{2n-2r}(X\times \mathcal{X},\mathbb{Q})$, and $\alpha_2$ is the restriction to $\mathcal{X}\times_B\mathcal{X}$ of a
class $\beta_2\in H^{2n-2r}( \mathcal{X}\times X, \mathbb{Q})$.

More precisely we can take $\beta_1\in \oplus_{i<n-r}H^{i}(X,\mathbb{Q})\otimes L^1H^{2n-2r-i}(\mathcal{X},\mathbb{Q})$,
and $\beta_2\in  \oplus_{i<n-r}L^1H^{2n-2r-i}(\mathcal{X},\mathbb{Q})\otimes H^{i}(X,\mathbb{Q})$, where $L$ stands for the Leray filtration on
$H^*(\mathcal{X},\mathbb{Q})$ with respect to the morphism $\pi:\mathcal{X}\rightarrow B$.
\end{lemm}
{\bf Proof.} Consider the smooth proper morphism
$$(\pi,\pi):\mathcal{X}\times_B\mathcal{X}\rightarrow B.$$
The relative  K\"{u}nneth decomposition gives
$$R^k(\pi,\pi)_*\mathbb{Q}=\bigoplus _{i+j=k} H^i_\mathbb{Q}\otimes H^j_\mathbb{Q},$$
where $H^i_\mathbb{Q}:=R^i\pi_*\mathbb{Q}$.
 The Leray spectral sequence of $(\pi,\pi)$, which  degenerates at $E_2$ (cf. \cite{delignedec}), gives the Leray filtration
  $L$ on $H^{2n-2r}(\mathcal{X}\times_B\mathcal{X},\mathbb{Q})$
 with graded pieces
 $$ Gr_L^lH^{2n-2r}(\mathcal{X}\times_B\mathcal{X},\mathbb{Q})=H^l(B,R^{2n-2r-l}(\pi,\pi)_*\mathbb{Q})
 =\bigoplus_{i+j=2n-2r-l}H^l(B,H^i_\mathbb{Q}\otimes H^j_\mathbb{Q}).$$
 Our assumption on $\alpha$ exactly says that it vanishes in the first quotient
 $H^0(B,R^{2n-2r}(\pi,\pi)_*\mathbb{Q})$, or equivalently, $\alpha\in L^1H^{2n-2r}(\mathcal{X}\times_B\mathcal{X},\mathbb{Q})$.
  Consider now the other graded pieces
 $$H^l(B,H^i_\mathbb{Q}\otimes H^j_\mathbb{Q}),\,\,{l>0,\,\,i+j=2n-2r-l}.$$
 Since $l>0$, and  $i+j=2n-2r-l$, we have  either $i<n-r$ or $j<n-r$.
 Let us consider the case where $i<n-r$: then the Lefschetz hyperplane section theorem tells that the sheaf
 $H^i_\mathbb{Q}$ is the constant sheaf with stalk $H^i(X,\mathbb{Q})$.
 Thus we find that $H^l(B,H^i_\mathbb{Q}\otimes H^j_\mathbb{Q})=H^i(X,\mathbb{Q})\otimes H^l(B,H^j_\mathbb{Q})$, which is
 a  Leray graded piece of $H^i(X,\mathbb{Q})\otimes H^{l+j}(\mathcal{X})$.
 Analyzing similarly the case where $j<n-r$, we conclude that
 the natural map
 $$ \bigoplus_{i<n-r}H^i(X,\mathbb{Q})\otimes L^1H^{2n-2r-i}(\mathcal{X},\mathbb{Q})\oplus \bigoplus_{j<n-r}L^1H^{2n-2r-j}(\mathcal{X},\mathbb{Q})\otimes H^{j}(X,\mathbb{Q})$$
 $$\rightarrow L^1H^{2n-2r}(\mathcal{X}\times_B\mathcal{X},\mathbb{Q})$$
 is surjective.
This proves the existence of the classes $\beta_1,\,\beta_2$.

\cqfd
In the case where $X$ has trivial Chow groups, we get an extra information:
\begin{lemm} \label{lebetaalgebraic}
With the same notations as above, assume that $X$ has trivial Chow groups
and that  $\alpha$ is the class of an algebraic cycle on $\mathcal{X}\times_B\mathcal{X}$. Then
we can choose  the $\beta_i$'s to be the restriction of   classes  of  algebraic cycles on
$B\times X\times X$.
\end{lemm}
{\bf Proof.} It suffices to show that we can choose the $\beta_i$'s to be classes of algebraic cycles
on $X\times \mathcal{X}$. Indeed,
$\mathcal{X}$ is a Zariski open set
in the natural fibration $$f:\mathbb{P}\rightarrow X,\,\mathbb{P}\subset \prod_i\mathbb{P}_i\times X,$$
$$\mathbb{P}:=\{(\sigma_1,\ldots,\sigma_r, x),\,\sigma_i(x)=0,\forall 1\leq i\leq r\}.$$
This is a fibration into  products of projective spaces, because we assumed the $L_i$'s are globally generated.
It follows that  $X\times \mathcal{X}$ is as well a Zariski open
set in  the corresponding fibration $X\times \mathbb{P}\rightarrow X\times X$
into  products of projective spaces.
The restriction map
$$ CH(X\times X\times \prod_{i}\mathbb{P}_i)\rightarrow CH(X\times \mathbb{P})$$
is then surjective, by the computation of the Chow groups of a projective bundle fibration (\cite[II,9.3.2]{voisinbook})
and thus, composing with the restriction  to the Zariski open set $X\times \mathcal{X}$, we get
$$ CH(X\times \mathbb{P})\rightarrow CH(X\times \mathcal{X})$$
is also surjective.
Hence the composition $CH(X\times X\times \prod_{i}\mathbb{P}_i)\rightarrow CH(X\times \mathcal{X})$, and a fortiori
the restriction map $CH(X\times X\times B)\rightarrow CH(X\times \mathcal{X})$ are surjective.

It remains to show that if $\alpha$ is algebraic, we can choose the $\beta_i$'s to be the restrictions of   classes  of  algebraic cycles on
$X\times \mathcal{X}$.

We have by the proof of Lemma
\ref{leleray} \begin{eqnarray}\label{equation11h38}
\alpha=\beta_{1\mid \mathcal{X}\times_B\mathcal{X}}+\beta_{2\mid \mathcal{X}\times_B\mathcal{X}},
\end{eqnarray}
where $\beta_1\in H^{*<n-r}(X,\mathbb{Q})\otimes L^1H^*(\mathcal{X},\mathbb{Q})$ and $\beta_2\in L^1H^*(\mathcal{X},\mathbb{Q})\otimes H^{*<n-r}(X,\mathbb{Q}) $.
We know that the cohomology of $X$ is generated by classes of algebraic cycles
$[z_{i,j}]\in H^{2i}(X,\mathbb{Q})$. Let us choose a basis $[z_{i,j}],\,2i<n-r$ of $H^{*<n-r}(X,\mathbb{Q})$.
Then we can choose cycle classes  $[z_{i,j}]^*\in H^{2n-2r-2i}(X,\mathbb{Q})$ in such a way that
the restricted classes $[z_{i,j}]^*_{\mid X_b}$ form the dual basis of $H^{*>n-r}(X_b,\mathbb{Q})$.
We can write by K\"{u}nneth's decomposition
\begin{eqnarray}
\label{autreeqn}
\beta_1=\sum_{i,j}[z_{i,j}]\cup \gamma_{i,j},\,\beta_2=\sum_{i,j}\gamma'_{i,j}\cup[z_{i,j}],\end{eqnarray}
where $\gamma_{i,j},\,\gamma'_{i,j}\in L^1H^* (\mathcal{X},\mathbb{Q})$.
Recalling that $\mathcal{X}\times_B \mathcal{X}\subset B\times X\times X$,
we denote by $p_{1,X}$ the  projection $\mathcal{X}\times_B \mathcal{X}\rightarrow  X$ to the first $X$ summand
and
$p_{2,X}$ the  projection $\mathcal{X}\times_B \mathcal{X}\rightarrow  X$
to the second $X$ summand.
Let $\pi'_2:\mathcal{X}\times_B \mathcal{X}\rightarrow \mathcal{X}$ be the second projection. Hence $\pi'_2$ is
a smooth projective morphism with fiber $X_b$  over any point of  $X_b\subset \mathcal{X}$.
The class $\alpha$ being algebraic, so is the class
$\pi'_{2*}(p_{1,X}^*[z_{i,j}]^*\cup \alpha)$. However, using (\ref{equation11h38}) and (\ref{autreeqn}), we have the equality
$$\gamma_{i,j}=\pi'_{2*}(p_{1,X}^*[z_{i,j}]^*\cup \alpha).$$ Indeed, we have $$\pi'_{2*}(p_{1,X}^*[z_{i,j}]^*\cup \beta_2)=
\pi'_{2*}(p_{1,X}^*[z_{i,j}]^*\cup (\sum_{i,j}\gamma'_{i,j}\cup[z_{i,j}]))=0$$

because $\gamma'_{i,j}\in L^1H^* (\mathcal{X},\mathbb{Q})$,
and
$$\pi'_{2*}(p_{1,X}^*[z_{i,j}]^*\cup \beta_1)=
\pi'_{2*}(p_{1,X}^*[z_{i,j}]^*\cup (\sum_{i,j}[z_{i,j}]\cup \gamma_{i,j}))=\gamma_{i,j}.$$
Similarly, $\pi'_{1*}(p_{2,X}^*[z_{i,j}]^*\cup \alpha)=\gamma'_{i,j}$ is algebraic on
$\mathcal{X}$. Thus all the $\gamma_{i,j}$'s and $\gamma'_{i,j}$'s are algebraic and so are
$\beta_1,\,\beta_2$.

\cqfd
{\bf Proof of Theorem \ref{theobrutintro}.} We keep  notation (\ref{notation}) and assume now that
the vanishing cohomology $H^{n-r}(X_b,\mathbb{Q})_{van}$ is supported on a codimension $c$
closed algebraic subset $Y_b\subset X_b$ for any $b\in B$. Consider the corrected (or vanishing) diagonal $\mathcal{D}_{van}$, which is a codimension
$n-r$ cycle of  $\mathcal{X}\times_B\mathcal{X}$ with $\mathbb{Q}$-coefficients.

By Lemma  \ref{coropropcrux}, it follows that there exists
a codimension $c$ closed algebraic subset $\mathcal{Y}\subset \mathcal{X}$ and
a codimension $n-r$ cycle $\mathcal{Z}$ on $\mathcal{X}\times_B\mathcal{X}$ with $\mathbb{Q}$-coefficients,
which is supported on $\mathcal{Y}\times_B\mathcal{Y}$ and such that
$$ [\mathcal{Z}_b]=[\mathcal{D}_{van,b}]=[\Delta_{b,van}],\,\,\forall b\in B.$$
Thus the class $[\mathcal{Z}]-[\Delta_{b,van}]\in H^{2n-2r}(\mathcal{X}\times_B\mathcal{X},\mathbb{Q})$
vanishes on the fibers $X_b\times X_b$.

Using Lemmas \ref{leleray} and \ref{lebetaalgebraic}, we conclude that there is
a cycle
$\Gamma\in CH^{n-r}(B\times X\times X)_\mathbb{Q}$ such that
\begin{eqnarray}\label{lapresquefin} [\mathcal{Z}]=[\mathcal{D}_{van}]+[\Gamma_{\mid \mathcal{X}\times_B\mathcal{X}}]\,\,{\rm in}\,\,
H^{2n-2r}(\mathcal{X}\times_B\mathcal{X},\mathbb{Q}).
\end{eqnarray}

 \begin{lemm}\label{lecyclemap}
   If $X$ has trivial Chow groups (that is, satisfies property $\mathcal{P}$),
the cycle class map $$CH^*(\mathcal{X}\times_B\mathcal{X})_\mathbb{Q}\rightarrow H^{2*}(\mathcal{X}\times_B\mathcal{X},\mathbb{Q})$$
is injective.
\end{lemm}
{\bf Proof.} Consider the blow-up $\widetilde{X\times X}$ of $X\times X$ along the diagonal. Applying Proposition \ref{propproduit}
and Lemma \ref{leblowup}, $\widetilde{X\times X}$ has trivial Chow groups.
A point of  $\widetilde{X\times X}$ parameterizes a
couple  $(x,y)$ of points of
$X$, together with a subscheme $z$ of length $2$ of $X$, with associated cycle $x+y$.
We thus  have the following natural variety
$$Q=\{(\sigma_1,\ldots,\sigma_r,x,y,z),\,\sigma_i\in \mathbb{P}_i,\,{\sigma_i}_{\mid z})=0,\,\, \forall i=1,\ldots, r\}\subset \prod_i\mathbb{P}_i\times\widetilde{X\times X}.$$
As the $L_i$'s are assumed to be very ample, the map
$Q\rightarrow \widetilde{X\times X}$ is a fibration with fibre over $(x,y,z)\in \widetilde{X\times X}$
 a product of projective spaces $\mathbb{P}_{i,z}$ of codimension $2$ in $\mathbb{P}_i$. By lemma
 \ref{leproj}, $Q$ also has trivial Chow groups.
 Let $Q_0\subset Q$ be the inverse image of $B$
 under the projection $Q\rightarrow \prod_i\mathbb{P}_i$. Then
 $Q_0$ is Zariski open in $Q$, so by Lemma \ref{leU}, the cycle map is also injective on cycles
 of $Q_0$. Finally, $Q_0$ maps naturally to $\mathcal{X}\times_B\mathcal{X}$ via the map
 $$\prod_i\mathbb{P}_i\times\widetilde{X\times X}\rightarrow \prod_i\mathbb{P}_i\times {X\times X}.$$
 The morphism $Q_0\rightarrow \mathcal{X}\times_B\mathcal{X}$ being projective and dominant, we conclude by Lemma
 \ref{ledominant} that the cycle map is injective on cycles of $\mathcal{X}\times_B\mathcal{X}$.

\cqfd

The proof  is then finished as follows.
From the equality (\ref{lapresquefin}) of cohomology classes, we deduce by the above lemma the following equality
of cycles:
 \begin{eqnarray}\label{lafin} \mathcal{Z}=\mathcal{D}_{van}+\Gamma_{\mid \mathcal{X}\times_B\mathcal{X}}\,\,{\rm in}\,\,
CH^{n-r}(\mathcal{X}\times_B\mathcal{X})_\mathbb{Q}.
\end{eqnarray}

We now fix $b$ and restrict this equality to $X_b\times X_b$. Then we find
$$ \mathcal{Z}_b=\Delta_{b,van}+\Gamma'_{\mid {X}_b\times{X}_b}\,\,{\rm in}\,\,
CH^{n-r}({X}_b\times{X}_b)_\mathbb{Q},
$$
where $\Gamma'\in CH(X\times X)_\mathbb{Q}$ is the restriction of $\Gamma$ to $b\times X\times X$.

Recalling that $\Delta_{b,van}=\Delta_{X_b}+\Gamma''_{\mid X_b\times X_b}$ for some codimension $n-r$-cycle
with $\mathbb{Q}$-coefficients $\Gamma''$ on $X\times X$, we conclude that
\begin{eqnarray}\label{egadiagb}\Delta_{X_b}=\mathcal{Z}_b+\Gamma_{1\mid X_b\times X_b},
\end{eqnarray}
where $\Gamma_1\in CH^{n-r}(X\times X,\mathbb{Q})$ and the cycle $\mathcal{Z}_b$ is by construction
supported on $\mathcal{Y}_b\times \mathcal{Y}_b$, with $\mathcal{Y}_b\subset X_b$ of codimension $\geq c$ for
general $b$.

Let $z\in CH_i(X_b)_\mathbb{Q}$, with $i<c$.
Then $(\mathcal{Z}_b)_*z=0$ since we may find a cycle rationally equivalent
to  $z$ in $X_b$ and   disjoint from $\mathcal{Y}_b$.
Applying both sides of (\ref{egadiagb}) to $z$ thus gives :
$$ z=(\Gamma_{1\mid X_b\times X_b})_*z\,\,{\rm in}\,\,CH_i(X_b)_\mathbb{Q}.$$
But it is obvious   that
$$(\Gamma_{1\mid X_b\times X_b})_*: CH(X_b)_\mathbb{Q}\rightarrow CH(X_b)_\mathbb{Q}$$
factors through $j_{b*}:CH(X_b)_\mathbb{Q}\rightarrow CH(X)_\mathbb{Q}$.
Now,  if $z$ is homologous to $0$ on $X_b$, $j_{b*}(z)$ is homologous to $0$ on
$X$, and thus it is rationally equivalent to $0$ on $X$ because $X$ has trivial Chow groups.
Hence we proved that the cycle map with $\mathbb{Q}$-coefficients is injective
on $CH_i(X_b)_\mathbb{Q} $ for $i<c$, which concludes the proof of the theorem.
\cqfd
\section{Variants and further applications \label{sec3}}
\subsection{Complete intersections with group action}
Theorem \ref{theobrutintro} applies to general  complete intersections
in projective space, the relation  (\ref{computconiveau}) giving the Hodge coniveau (hence conjecturally the geometric
coniveau $c$).
There are interesting variants coming from the study of
complete intersections $X_b$ of $r$ hypersurfaces in projective space $\mathbb{P}^n$, or in a product of projective spaces, invariant under a finite group action.
Let $G$ acts on $X_b$, and let $\chi:G\rightarrow \mathbb{Z}/2\mathbb{Z}=\{1,-1\}$ be a character of
$G$. Then consider the sub-Hodge structure
$$L^\chi=\{\alpha\in H^{n-r}(X_b,\mathbb{Q})_{prim},\,g^*\alpha=\chi(g)\alpha,\,\forall g\in G\}\subset H^{n-r}(X_b,\mathbb{Q})_{prim}.$$

In general, it has a larger coniveau than $X_b$. For example if
$X_b$ is a quintic surface in $\mathbb{P}^3$, defined by an invariant polynomial under the linearized group action
of $G\cong \mathbb{Z}/5\mathbb{Z}$ with generator $g$
on $\mathbb{P}^3$ given by
$$g^*X_i=\zeta^iX_i,\,i=0,\ldots, 3,$$
where $\zeta $ is a nontrivial $5$-th root of unity,
then $H^2(S,\mathbb{Q})^{inv}$ has no $(2,0)$-part, hence is of coniveau $1$, while $H^{2,0}(S)\not=0$ so the coniveau of
$H^2(S,\mathbb{Q})_{prim}$ is $0$. The quotient surfaces $S/G$ is a quintic Godeaux surface (cf. \cite{voisinscuola}).

Note that the Hodge structure $L^\chi$ corresponds to the projector
$\frac{1}{|G|}\sum_{g\in G}\chi(g)g^*$ acting on $L$, and it is given by the action of the
$n-r$-cycle
$$\Gamma_\chi:=\sum_{g\in G}\chi(g)\Delta_{b,van,g},$$
where $\Delta_{b,van,g}=(Id,g)_*(\Delta_{b,van})\in CH_{n-r}(X_b\times X_b)_\mathbb{Q}$.
The generalized Bloch conjecture \ref{blochconj} (extended to  motives)
predicts the following :
\begin{conj}\label{congroup} Assume $L^\chi$ has coniveau $\geq c$. Then the cycle map is injective on $CH(X_b)_\mathbb{Q}^\chi$
for $i<c$.
\end{conj}
If $\chi$ is the trivial character, this conjecture is essentially equivalent to the previous one
by considering $X/G$ or a desingularization of it. Even in this case, one needs
to make assumptions on the linearized group action in order to apply
the same strategy as in the proof of Theorem \ref{theobrutintro}.
The case of non trivial character cannot be
reduced to the previous case.

In order to apply a strategy similar to the one applied
for the proof of Theorem \ref{theobrutintro}, we need some assumptions. Indeed, if the group $G$ is too big,
like the automorphisms group of the Fermat hypersurface, there are to few invariant complete intersections
to play on the geometry of the universal family $\mathcal{X}\rightarrow B$ of $G$-invariant
complete intersections.

In any case, what we get mimicking the proof of Theorem \ref{theobrutintro} is the following:
$X$ is as before a smooth projective variety of dimension $n$ satisfying property $\mathcal{P}$ and $G$ is a finite group acting
on $X$.
We study complete intersections $X_b\subset X$
of $r$ $G$-invariant ample hypersurfaces $X_i\in |L_i|^G$ : Let $G$ acts via the character $\chi_i$ on the considered
component of $|L_i|^G$. The basis $B$ parameterizing such
complete intersections is thus a Zariski open set in $\prod_i\mathbb{P}(H^0(X,L_i)^{\chi_i})$. As before
we denote by $\mathcal{X}\rightarrow B$ the universal complete intersection.
\begin{theo}\label{variantgoupaction} Assume

(i) The variety $\mathcal{X}\times_B\mathcal{X}$ satisfies property $\mathcal{P}_{n-r}$.

(ii) The Hodge structure on $H^{n-r}(X_b,\mathbb{Q})_{van}^{\chi}$ is supported on a closed algebraic subset $Y_b\subset X_b$
of codimension $c$. (Conjecturally, this is satisfied if the Hodge coniveau of $H^{n-r}(X_b,\mathbb{Q})_{van}^{\chi}$ is $\geq c$,
 cf. Conjecture \ref{hodgeconjgen}).

(iii) Conjecture \ref{conjstandard} holds for codimension $n-r$ cycles.

Then the cycle map $CH_i(X_b)^\chi_\mathbb{Q}\rightarrow H^{2n-2r-2i}(X_b,\mathbb{Q})^\chi$ is injective
for any $b\in B$.
\end{theo}

\begin{rema}\label{remasurfacesthreefolds}
{\rm In the case where $X_b$ are  surfaces with $h^{2,0}(X_b)^\chi=0$, the assumption (ii) is automatically satisfied, by the same arguments
as in the  proof of Theorem \ref{surfaces}. We thus get an alternative proof of
the main theorem  of \cite{voisinscuola}, where the Bloch conjecture is proved for the general Godeaux surfaces
(quotients of quintic surfaces by a free action of  $\mathbb{Z}/5\mathbb{Z}$, or quotients of complete intersections of four quadrics in $\mathbb{P}^6$
by a free action of  $\mathbb{Z}/8\mathbb{Z}$).

In the case of threefolds $X_b$ of Hodge coniveau $1$, we can also conclude
that $CH_0(X_b)^\chi_0=0$ if  (i) is satisfied and  the generalized Hodge conjecture  is satisfied by
the coniveau $1$ Hodge structure on $H^3(X_b,\mathbb{Q})^\chi$.
 Indeed,
we used  conjecture
\ref{conjstandard}  in the proof essentially for the proof of Lemma  \ref{leeasybutcrucial}, which says that
 if a certain Hodge structure $L\subset H^*(X_b,\mathbb{Q})$ is supported on a codimension $c$ closed algebraic subset
  $Y_b$, the corresponding projector
 has a class which comes from the class of a cycle supported in $Y_b\times Y_b$. This
 will be satisfied if ${\rm dim}\,X_b=3$, $L\subset H^3(X,\mathbb{Q})^\chi$ supported on $Y_b\times Y_b$
  because we know that the degree $6$ Hodge class of the projector $\pi_L$
 is supported on the codimension
 $2$ closed algebraic subset $Y\times Y$ (or rather a desingularization of it), so that we can apply Lemma
 \ref{lemmacodim2}.

 This way the second  result of \cite{voisinscuola} (quintic hypersurfaces with involutions)
 and the main application of \cite{peters} ($3$-dimensional complete intersection in weighted projective space) are reproved : in both cases
   we are reduced to prove the generalized Hodge conjecture
 for the coniveau $1$ Hodge structure on their cohomology of degree $3$.

}
\end{rema}
\begin{example} {\rm Consider
a Calabi-Yau hypersurface $X_f$ in $X=\mathbb{P}^n$  defined by an equation
$f$ invariant by the involution
$i:i^*(X_0,\ldots,X_n)=(-X_0,-X_1,X_2,\ldots, X_n)$. Then
$H^{n-1}(X)_{prim}^-$ has coniveau $1$, since
$i$ acts by $Id$ on $H^{n-1,0}(X_f)$.
In \cite{voisinscuola}, the case of $3$-dimensional quintics is studied, and it is proved there that
$i$ acts by $Id $ on $CH_0(X)$  in this case. One step is the proof that the generalized Hodge conjecture
holds for the coniveau $1$ Hodge structure  $H^3(X,\mathbb{Q})^-$. Having this, Theorem \ref{variantgoupaction} gives a drastically
simplified proof of this result.
}
\end{example}
\begin{example}{\rm The following class of examples is constructed in \cite{ctv}: $X=\mathbb{P}^1\times \mathbb{P}^3$,
with the following group action: $G=\mathbb{Z}/5\mathbb{Z}$ acts on $\mathbb{P}^1\times \mathbb{P}^3$ in the following linearized way:
Let $g$ be a generator of $G$ and $\zeta $  a nontrivial $5$-th root of unity. Then
if $x,y$ are homogeneous coordinates on
$\mathbb{P}^1$ and $x_0,x_1,x_2,x_3$ are homogeneous coordinates on
$\mathbb{P}^3$, we set:
$$g^*x=x,\,g^*y=\zeta y,$$
$$g^*x_i=\zeta^ix_i,\,i=0,\ldots, 3.$$

 We then consider hypersurfaces $X_f\subset \mathbb{P}^1\times \mathbb{P}^3$
 defined by an equation $f=0$ of  bidegree $(3,4)$, where
 $f\in H^0(\mathbb{P}^1\times \mathbb{P}^3,\mathcal{O}_{\mathbb{P}^1\times \mathbb{P}^3}(3,4))^G$.

These hypersurfaces $X_f$ have a few ordinary quadratic singularities. The varieties
$X'_f$ obtained as a desingularization of $X_f/G$ have  $h^{3,0}(X'_f)=0$ (and also $h^{i,0}(X'_f)=0$ for
$i=1,2$).
For the general such  variety, Theorem \ref{variantgoupaction} tells that the generalized Hodge conjecture
for $H^3(X'_f,\mathbb{Q})$ implies (and in fact is equivalent to by Theorem \ref{theotoutlemonde})
the equality $CH_0(X'_f)=\mathbb{Z}$.
The interest in these examples comes from the fact, proved in \cite{ctv}, that the
Hodge conjecture is not satisfied for  integral Hodge classes of degree $4$ on $X'_f$.

}
\end{example}
\subsection{Self-products}
Let $Y$ be a smooth projective variety. We will assume for simplicity
that $H^{i,0}(Y)=0$ for
$i\not=0,\,m:={\rm dim}\,Y$. (This will be the case if $Y$ is a complete intersection
of ample hypersurfaces in a projective variety with trivial Chow groups.)

 \begin{lemm}\label{leconiveauCY} For $k>p_g(Y)=h^{m,0}(Y)$, the Hodge structure of
 weight $km$ on $\bigwedge^kH^m(Y,\mathbb{Q} )$
 has coniveau $\geq1$. In particular, if $h^{m,0}(Y)=1$, the Hodge structure of
 weight $2m$ on $\bigwedge^2H^m(Y,\mathbb{Q} )$ has coniveau $\geq1$.
 \end{lemm}
 {\bf Proof.} Indeed, the $(km,0)$-piece of the Hodge structure on
$ \bigwedge^kH^m(Y,\mathbb{Q} )$ is equal to
$\bigwedge^{k}H^{m,0}(Y)$, hence it is $0$ for $k>h^{m,0}(Y)$.

 \cqfd
  Conjecture \ref{blochconj} (or rather its generalization to motives)
 predicts the following (see below for more detail):
 \begin{conj} Assume $Y$ satisfies the above assumption and has $h^{m,0}=1$.
 Then, for any $z,\,z'\in CH_0(Y)$ with ${\rm deg}\,z={\rm deg}\,z'=0$, one has
 $z\times z'-z'\times z=0$ in $CH_0(Y\times Y)$ for $m$ even and
 $z\times z'+z' \times z=0$ in $CH_0(Y\times Y)$ for $m$ odd
 \end{conj}
 The case $m=2$ is particularly interesting, as noticed in \cite{voisinsymmetric}. In this case, we indeed have :

 \begin{lemm}\label{leK3} Let $H, H^{p,q}$ be a  weight $2$ Hodge structure of  $K3$ type, namely $h^{2,0}=1$.
 Then the Hodge structures on $\bigwedge ^{2k}H$ all have niveau $\leq 2$ (that is coniveau $\geq k-1$).

 \end{lemm}
 {\bf Proof.} Write $H=H^{2,0}\bigoplus H^{1,1}\oplus H^{0,2}$.
 Then $$\bigwedge^kH=H^{2,0}\otimes\bigwedge^{k-1}H^{1,1}\oplus (\bigwedge^{k}H^{1,1}\oplus H^{2,0}\otimes H^{0,2}\otimes \bigwedge^{k-2}H^{1,1})\oplus \bigwedge^{k-1}H^{1,1}\otimes H^{0,2}$$
 is the Hodge decomposition of $\bigwedge^kH$, whose first nonzero term is of type $(k+1,k-1)$.
 \cqfd
When $k>{\rm dim}\,H$, we of course have that the Hodge structure on $\bigwedge^kH$ is trivial.
Applying these observations to the case where $H=H^2(S,\mathbb{Q})$ where
$S$ is  an algebraic  $K3$ surface, we find that
Conjecture \ref{blochconj} (or rather, its extension to motives)
predicts the following (cf. \cite{voisinsymmetric}):
\begin{conj} \label{conjproductk3} (i) Let $S$ be an algebraic $K3$ surface. Then for any
$k\geq 2$, and $i\leq k-2$,  the projector $\pi_{alt}=\sum_{\sigma\in \mathfrak{S}_k}(-1)^{\epsilon(\sigma)}
\sigma \in CH^{2k}(S^k\times S^k)$
composed with the Chow-K\"unneth projector $\pi_{2}^{\otimes k}$ (cf. \cite{murre})
acts as  $0$ on $CH_{i}(S^k)_\mathbb{Q}$ for $i\leq k-2$.

(ii) For $k>b_2(S)$, this projector is identically $0$.
\end{conj}

Note that (ii) above is essentially Kimura's finite dimensionality conjecture
\cite{kimura} and applies to any regular surface. One may wonder whether it could be attacked by
the methods of the present paper for the case of quartic $K3$ surfaces. The question would be
essentially to study whether the fibered product
$\mathcal{X}^{2k/B}$ of the universal such $K3$ surface satisfies property
$\mathcal{P}$.
For small $k$ this is easy, but we would need to know this in the range $k\geq 44$ in order to prove the Kimura conjecture.
This seems to be very hard.

The fact that this is true for small $k$ (see below) shows that Conjecture \ref{conjproductk3} is implied
by  the generalized Hodge conjecture for the self-products $S^k$ and the
coniveau $k-1$ Hodge structures
$\bigwedge^kH^2(S,\mathbb{Q})\subset H^{2k}(S^k,\mathbb{Q})$.

Let us be a little more explicit  in  the case of general Calabi-Yau complete intersections  and for $k=2$.
Let $X_b$ be a smooth Calabi-Yau complete intersection of dimension $m$  in projective space $\mathbb{P}^n$.
Let $\Delta_{b,van}\in CH^{m}( X_b\times X_b)_\mathbb{Q}$ be the corrected diagonal,
whose action on $H^*(X_b,\mathbb{Q})$ is the projection on $H^{m}(X_b,\mathbb{Q})_{van}$.
On $X_b\times X_b\times X_b\times X_b$,
there is the induced  $2m$-cycle
$$\Delta_{b,van,2}:= p_{13}^*\Delta_{b,van}\cdot p_{24}^*\Delta_{b,van},$$
where $p_{ij}$ is the projection from $X_b^4$ to the product $X_b^2$ of its $i$-th and $j$-th factor.
 The  action on $\Delta_{b,van,2}$   seen as a self-correspondence of $X_b^2$  on $H^*(X_b^2,\mathbb{Q})$ is the orthogonal projector on
 $$p_1^* H^{m}(X_b,\mathbb{Q})_{van}\otimes p_2^*H^{m}(X_b,\mathbb{Q})_{van}\subset H^{2m}(X_b\times X_b,\mathbb{Q}).$$
 If instead of $\Delta_{b,van,2}$, we consider
 $$\Delta_{b,van,2}^\tau:= p_{14}^*\Delta_{b,van}\cdot p_{23}^*\Delta_{b,van},$$
 then the action on $\Delta_{b,van,2}$   seen as a self-correspondence of $X_b^2$  on $H^*(X_b^2,\mathbb{Q})$ is the composition
 of the previous projector with the permutation
 $$\tau_*: H^{m}(X_b,\mathbb{Q})_{van}\otimes H^{m}(X_b,\mathbb{Q})_{van}\rightarrow H^{m}(X_b,\mathbb{Q})_{van}\otimes H^{m}(X_b,\mathbb{Q})_{van} $$
 exchanging summands. Note that the inclusion
 $$H^{m}(X_b,\mathbb{Q})_{van}\otimes H^m(X_b,\mathbb{Q})_{van}\subset H^{2m}(X_b\times X_b,\mathbb{Q})$$
 sends the antiinvariant part on the left to the antiinvariant part under $\tau$ on the right if $m$ is even, and
 to the invariant part under $\tau$ on the right if $m$ is even. This is due to the fact that the cup-product
 on cohomology is graded commutative.

 Hence we conclude that
 $$\Delta_{b,van,2}^\sharp:={\Delta_{b,van,2}-\Delta_{b,van,2}^\tau}$$
 acts on $H^*(X_b^2,\mathbb{Q})$ as twice the projector onto $\bigwedge^2H^m(X_b\times X_b,\mathbb{Q})_{van}$ if $m$ is even,
 and that
 $$\Delta_{b,van,2}^{inv}:={\Delta_{b,van,2}+\Delta_{b,van,2}^\tau}$$
acts on $H^*(X_b^2,\mathbb{Q})$ as twice the projector onto $\bigwedge^2H^m(X_b\times X_b,\mathbb{Q})_{van}$ if $m$ is odd.

In both cases, using Lemma \ref{leconiveauCY}, we get that this is twice the orthogonal  projector associated to a sub-Hodge structure
of coniveau $\geq1$.

Restricting to the case of Calabi-Yau hypersurfaces in $\mathbb{P}^n$, (so $m=n-1$),  an easy adaptation of
the proof of Theorem \ref{theobrutintro} gives now:
\begin{theo} (cf. Theorem \ref{theoproduitCY})
Assume Conjecture \ref{conjstandard} and the generalized conjecture holds for the coniveau $1$ Hodge structure on
$\bigwedge^2H^{n-1}(X_b\times X_b,\mathbb{Q})_{van}\subset H^{2n-2}(X_b\times X_b,\mathbb{Q})$, where $X_b$ is a general Calabi-Yau
hypersurface in projective space.
Then the general such $X_b$ has the following property:

(i) If $n-1$ is even, for any two $0$-cycle $z,\,z'$ of degree $0$ on $X_b$, we have $z\times z'-z'\times z=0$ in $CH_0(X\times X)$.

(ii) If $n-1$ is odd, for any two $0$-cycle $z,\,z'$ of degree $0$ on $X_b$, we have $z\times z'+z'\times z=0$ in $CH_0(X\times X)$.
\end{theo}
{\bf Proof.} We just sketch the proof, as it is actually a variant of the proof of Theorem \ref{theobrutintro}.
With the same notations as in \ref{notation} (where $X$ will be the projective space $\mathbb{P}^{n}$),
we claim that it suffices  to show that, if $n-1$ is even,
the spread-out cycles
$$\mathcal{D}_{van,2}^\sharp:={\mathcal{D}_{van,2}-\mathcal{D}_{van,2}^\tau}\in
CH^{2n-2}(\mathcal{X}\times_B\mathcal{X}\times_B\mathcal{X}\times_B\mathcal{X})_\mathbb{Q}$$
can be written as a sum
\begin{eqnarray}\label{eqndimanche0} \mathcal{D}_{van,2}^\sharp=
\mathcal{Z}_1+\mathcal{Z}_2\,\,{\rm in}\,\,CH^{2m}(\mathcal{X}\times_B\mathcal{X}\times_B\mathcal{X}\times_B\mathcal{X})_\mathbb{Q},
\end{eqnarray}
where $\mathcal{Z}_1$ is supported on $\mathcal{Y}\times_B\mathcal{Y}$, with $\mathcal{Y}\varsubsetneqq \mathcal{X}\times_B\mathcal{X}$,
and $\mathcal{Z}_2$ is a cycle which is the restriction of cycles
on various copies  of $\mathcal{X}\times_B\mathcal{X}\times B\times X\times X$ (ordered adequately), via the
inclusion $ \mathcal{X} \subset B\times X$; similarly for $n-1$ odd,
with $\mathcal{D}_{van,2}^{\sharp}$ replaced  by
 $$\mathcal{D}_{van,2}^{inv}:={\mathcal{D}_{van,2}+\mathcal{D}_{van,2}^\tau}\in
CH^{2n-2}(\mathcal{X}\times_B\mathcal{X}\times_B\mathcal{X}\times_B\mathcal{X})_\mathbb{Q}.$$

Indeed, if we know this, restricting to a general point $b\in B$, we get that for $n-1$ even
$\Delta_{van,2}^\sharp:={\Delta_{van,2}-\Delta_{van,2}^\tau}\in CH^{2n-2}(X_b^4)_\mathbb{Q}$ can be written as a sum
\begin{eqnarray}\label{eqndimanche} \Delta_{van,2}^\sharp={Z}_1+{Z}_2,
\end{eqnarray}
where ${Z}_1$ is supported on $\mathcal{Y}_b\times\mathcal{Y}_b$, with $\mathcal{Y}_b\varsubsetneqq {X}_b\times {X}_b$,
and ${Z}_2$ is a cycle which is the restriction of cycles
on various copies  of ${X}_b\times{X}_b\times X\times X$ (ordered adequately), via the
inclusion $j_b:{X}_b \hookrightarrow X=\mathbb{P}^{n}$. Similarly for $n-1$ odd with $\Delta_{van,2}^\sharp$ replaced by
$\Delta_{van,2}^{inv}$.
We see equation (\ref{eqndimanche})
as an equality of self-correspondences of
$X_b^2$ and we let  both sides of (\ref{eqndimanche}) act on $z\times z'$, where $z,\,z'\in CH_0(X_b)$ have degree $0$.
On the left, we get $\Delta_{van,2}^\sharp(z\times z')=z\times z'-z'\times z$. Next the cycle
${Z}_1$ being supported on $\mathcal{Y}_b\times\mathcal{Y}_b$, with $\mathcal{Y}_b\varsubsetneqq {X}_b\times {X}_b$, acts
trivially on $CH_0(X_b\times X_b)$.
We thus get (for $n-1$ even), decomposing $Z_2$ as a sum $Z_2=\sum_i Z_{2,i\mid X_b^4}$ where
$Z_{2,i}\in CH^{2n-2}(X_b\times \ldots \times X\times\ldots X_b)$,
with the factor $X$  put in $i$-th position:
\begin{eqnarray}z\times z'-z'\times z=\sum_i(Z_{2,i\mid X_b^4})_*(z\times z')\,\,{\rm in}\,\,CH_0(X^2)_\mathbb{Q}.
\end{eqnarray}

If $i=1,\,2$, $(Z_{2,i\mid X_b^4})_*(z\times z')$ vanishes because, as $CH_0(X)=CH_0(\mathbb{P}^{n})=\mathbb{Z}$,
both cycles $z\times j_{b*}z'\in CH_0(X_b\times X)$ and $j_{b*}z\times z'\in CH_0(X\times X_b)$
vanish.
For $i=3,\,4$, we have  that $(Z_{2,i\mid X_b^4})_*(z\times z')$ belongs to $CH_1(X_b\times X)_{\mid X_b\times X_b}$ or to
$CH_1(X\times X_b)_{\mid X_b\times X_b}$ with $X=\mathbb{P}^{n}$. Thus we get using the decomposition
of $CH_1(X_b\times \mathbb{P}^{n})$ as $ CH_0(X_b)\otimes \mathbb{Z}h_1\oplus CH_1(X_b)\times \mathbb{Z}h_0$, where $h_1$ is the class of a
line and $h_0$ is the class of a point  in $\mathbb{P}^{n}$,  an equality:
$$ z\times z'-z'\times z=w_1\times h_{1\mid X_b}+h_{1\mid X_b}\times w_2\,\,{\rm in}\,\,CH_0(X_b\times X_b)_\mathbb{Q},$$
with $w_i\in CH_0(X_b)$.
Applying $pr_{1*}$ and $pr_{2*}$ to both sides of this equality, we finally get  that $w_1=w_2=0$ in $CH_0(X_b)_\mathbb{Q}$.
Thus we proved  assuming (\ref{eqndimanche0})
that $z\times z'-z'\times z=0$ in $CH_0(X^2_b)_\mathbb{Q}$ for $n-1$ even and the same proof will show that
$z\times z'+z'\times z=0$ in $CH_0(X^2_b)_\mathbb{Q}$ for $n-1$ odd. As these cycles belong to the kernel of the Albanese map,
we also conclude by Roitman's theorem \cite{roitman} that these equalities in fact hold in $CH_0(X^2_b)$.

It remains to see how to get (\ref{eqndimanche0}) from the condition that the generalized Hodge
conjecture holds for the coniveau $1$ Hodge structure on $\bigwedge^2H^{n-1}(X_b)_{prim}$ combined
with  Conjecture \ref{conjstandard}. As in the proof of Theorem \ref{theobrutintro}, we find that under these two assumptions,
we have an equality of cycle classes
\begin{eqnarray}\label{eqndimanche0coh}[\mathcal{D}_{van,2}^\sharp]= [\mathcal{Z}_1]+[\mathcal{Z}_2]\,\,{\rm in}\,\,H^{4n-4}(\mathcal{X}\times_B\mathcal{X}\times_B\mathcal{X}\times_B\mathcal{X},\mathbb{Q}),
\end{eqnarray}
where $\mathcal{Z}_1$ is supported on $\mathcal{Y}\times_B\mathcal{Y}$, with $\mathcal{Y}\varsubsetneqq \mathcal{X}\times_B\mathcal{X}$,
and $\mathcal{Z}_2$ is a cycle which is the restriction of cycles
on various copies  of $\mathcal{X}\times_B\mathcal{X}\times B\times X\times X$.

Equation (\ref{eqndimanche0}) follows from (\ref{eqndimanche0coh}) and from the following Proposition \ref{nextprop}.
This finishes the proof of Theorem \ref{theoproduitCY}.
\cqfd
\begin{prop} \label{nextprop}The fourth fibered product $\mathcal{X}\times_B\mathcal{X}\times_B\mathcal{X}\times_B\mathcal{X}$
of universal hypersurfaces of degree $\geq 3$ in $\mathbb{P}^n$ satisfies property $\mathcal{P}_{2n-2}$.

\end{prop}
{\bf Proof.}  As we are interested into cycles of codimension $\leq 2n-2$, we can restrict to the open
set $\mathcal{X}^{4/B}_0$ defined as the complement of the small relative diagonal $\mathcal{X}\subset \mathcal{X}^{4/B}_0$
which is of codimension $3n-3$.

We apply the relative version of
Lemma \ref{ledesing} below. This provides us with a certain   blow-up  $\widetilde{\mathcal{X}^{4/B}_0}$ of the
relative  Fulton-MacPherson
configuration space  (cf. \cite{fulmph}). It is   smooth and proper over
${\mathcal{X}^{4/B}_0}$. In order to prove the result, it suffices by Lemma \ref{ledominant}
to show that $\widetilde{\mathcal{X}^{4/B}_0}$
satisfies property $\mathcal{P}_{2n-2}$. By the functoriality statement in   Lemma
\ref{ledesing}, there is a morphism
from $\widetilde{\mathcal{X}^{4/B}_0}$ to $\widetilde{({\mathbb{P}^n})^{4}_0}$ and in particular
to the  punctual Hilbert scheme
${\rm Hilb}^4(\mathbb{P}^n)$, so that an element of $\widetilde{\mathcal{X}^{4/B}_0}$  determines
a $4$-uple $(x_1,\ldots,x_4)\in X_b$ together with a subscheme $z$ of $X_b\subset \mathbb{P}^n$ of length $4$ with associated cycle
$x_1+\ldots+x_4$.
It is an easy result that any subscheme $z$ of length $4$ of $\mathbb{P}^n$ whose support
consists of at least two points imposes independent conditions to degree $n+1$ hypersurfaces, with $n+1\geq 3$.
It follows that
$\widetilde{\mathcal{X}^{4/B}_0}$ can be realized as a Zariski open set of
a projective bundle over the space
$\widetilde{({\mathbb{P}^n})^{4}_0}$
constructed below. Namely, over a point $u$ in this  space, giving rise to
$x_1,\dots, x_4$ together with a schematic structure $z$ with associated cycle $x_1+\ldots+x_4$, the fiber is the
projective space $\mathbb{P}(H^0(\mathbb{P}^n,\mathcal{I}_z(n+1)))$, and
we have to take inside it the Zariski open set
which parameterizes smooth hypersurfaces.

Lemma \ref{leblowup} says  that $\widetilde{({\mathbb{P}^n})^{4}_0}$ satisfies property $\mathcal{P}_{2n-2}$. By Lemma \ref{leproj}, it follows
that the  projective bundle described above over $\widetilde{({\mathbb{P}^n})^{4}_0}$  also satisfies
$\mathcal{P}_{2n-2}$. By Lemma \ref{leU}, the Zariski open set $\widetilde{\mathcal{X}^{4/B}_0}$ inside it also
satisfies
$\mathcal{P}_{2n-2}$.
\cqfd
\begin{lemm}\label{ledesing} (cf. \cite{lehnsorger})
Let $X$ be a smooth variety of any dimension
$n$. Denote by $X^4_0$ the open set $X^4\setminus \Delta_{X,4}$, where $\Delta_{X,4}\cong X$
is the small diagonal. There is a smooth variety $\widetilde{{X}^{4}_0}$
which admits a morphism to ${\rm Hilb}^4({X})$ (so the rational
map  $\sigma:{{X}^{4}_0}\dashrightarrow{\rm Hilb}^4({X})$ is desingularized on $\widetilde{{X}^{4}_0}$),
whose construction is  functorial under  immersions, and which satisfies property $\mathcal{P}_{2n-2}$
if $X$ satisfies property $\mathcal{P}$.
\end{lemm}
{\bf Proof.} We will just describe  the construction of $\widetilde{{X}^{4}_0}$ over  a neighborhood in $X^4_0$ of a point
of ${X}^{4}_0$ of type $(3,1)$, that is a point
which corresponds to a relative $0$-cycle of the form $3x+y$. The case of
points of type $(2,2)$ (that is a point
which corresponds to a  $0$-cycle of the form $2x+2y$) or with a support of cardinal $\geq 3$
is easy  and left to the reader. As the set of points of type $(3,1)$ and $(2,2)$ are disjoint in $X^4_0$ there is no problem
to glue the local constructions.

  We first blow-up inside ${X}^{4}_0$
the union of the images under  permutation
of the diagonals of type  $\Delta_{xxxy}$ parameterizing
  the points
 $$(x,x,x,y)\in {X}^{4}_0,\,x,\,y\in {X},\,x\not=y.$$
 We then blow-up the (disjoint union of the) proper transforms of the images under permutations
 of the  big diagonal $\Delta_{xxyz}$
which is defined as
  the closure of the set of points
 $$(x,x,y,z)\in {X}^{4}_0,\,x,\,y,\,z\in {X},\,x,\,y,\,z\,\,{\rm distinct}.$$
 What we get at this point is nothing but the Fulton-MacPherson compactification $X(4)_0$
 of
 the configuration space of $4$ points, at least over the  open set $X^4_0$
 of $X^4$ (cf. \cite{fulmph}).

 The rational map $\sigma:X(4)_0\dashrightarrow {\rm Hilb}^4(X)$ is not yet a morphism, as shown to us
 by Totaro. What we need to blow-up is the following locus $M$ (pointed out by Totaro):
 over the  diagonal $\Delta_{x,x,x,y}$ (or any image  of it under permutation),
 the exceptional divisor over $\Delta_{x,x,x,y}$ is (before the second blow-up)
 the projective bundle $\mathbb{P}(p_x^*T_X\oplus p_x^*T_X)$ where $p_x$ is the first
 projection
 $\Delta_{xxxy}\cong X\times X\rightarrow X$.
 Let $M'\subset \mathbb{P}(p_x^*T_X\oplus p_x^*T_X)$ be the locus of couples $(v_1,v_2)$ where $v_1,\,v_2$ are colinear.
 $M'$ is  isomorphic to $\mathbb{P}^1\times  \mathbb{P}(p_x^*T_X)$.
 Let $M$ be the proper transform of $M'$ under the second blow-up, that is in  $X(4)_0$. $M$ is isomorphic
 to $M'$.

It is explained in the letter \cite{lehnsorger} that $\sigma$ becomes well-defined
 on the blow-up of $X(4)_0$ along $M$, as a consequence of Hayman's theorem \cite{hay}. This of course concerns the case where
 ${\rm dim}\,X=2$. However,  we are looking at ${\rm Hilb}^4(X)_0$, the open set of ${\rm Hilb}^4(X)$
 where the support has cardinality at least $2$, and in fact are mostly concerned with
 the neighborhood in ${\rm Hilb}^4(X)$ of punctual subschemes  of type $(3,1)$. As this is locally
 (in the \'{e}tale or analytic topology)  isomorphic to
 the product $X\times {\rm Hilb}^3(X)$, we are reduced to study
 the  case of ${\rm Hilb}^3(X)$ in a neighborhood of a fat point $z$.
 We observe now that
 any length $3$ subscheme $z\subset X$ is contained in a smooth surface
 in $X$. More precisely, if
  $n$ is the dimension of $X$, we choose a linear system of hypersurfaces
 $H_0\ldots,H_n$ in $X$ with the property that $z$ imposes $3$ independent
 conditions to $<H_0,\ldots,H_n>$, and that the locus $\Sigma_z\subset X$ defined
 by the linear subsystem $I_z\subset <H_0,\ldots,H_n>$ is smooth. The map
 $$\phi: {\rm Hilb}^3(X)\rightarrow {\rm Grass} (3, n+1),$$
 $$z'\mapsto I_{z'} \subset <H_0,\ldots,H_n>$$
 is then well-defined near $z$, dominant, with fibers through a point  $z'\in {\rm Hilb}^3(X)$ close to $z$ the Hilbert scheme
 ${\rm Hilb}^3(\Sigma_{z'})$ of the smooth surface
 $\Sigma_{z'}$. Using this, we easily reduce the general case to the surface case.

 Applying Lemmas \ref{leblowup} and \ref{leU}, it is easy to show that the resulting variety
 $\widetilde{{X}^{4}_0}$,
 satisfies property $\mathcal{P}_{2n-2}$
if $X$ satisfies property $\mathcal{P}$.
\cqfd

\end{document}